\documentclass[journal]{IEEEtran}
\IEEEoverridecommandlockouts
\usepackage{lineno}
\usepackage[cmex10]{amsmath}
\usepackage{amscd,amsfonts,amsmath,amssymb,mathrsfs,pifont,stmaryrd,tipa}
\usepackage{array,cases,dsfont,graphicx,texdraw}
\usepackage{graphics}
\usepackage{epstopdf}
\usepackage{epsfig}
\usepackage{cite}
\usepackage{mathbbold}
\usepackage{stfloats}
\usepackage{subfig}
\usepackage{wrapfig}
\usepackage{amsfonts}
\usepackage{amssymb}
\usepackage{graphicx}%
\usepackage{xcolor}
\usepackage{caption}
\usepackage[labelfont={color=black}]{caption}

\newtheorem{theorem}{\bf{Theorem}}

\newtheorem{assumption}{Assumption}

\newtheorem{proposition}{Proposition}
\newenvironment{proof}{\noindent{\em Proof:\/}}{\hfill $\Box$\par}
\newtheorem{remark}{\bf{Remark}}

\allowdisplaybreaks
\begin{document}

\title{Distributed Stochastic Optimization under Heavy-Tailed Noises}
\author{Chao Sun, Huiming Zhang, Bo Chen, Li Yu\thanks{Chao Sun, Bo Chen and Li Yu are with the  Department of Automation, Zhejiang University of Technology, China. Huiming Zhang is with the Institute of Artificial Intelligence, Beihang University, China. Email: bchen@aliyun.com.}}
\maketitle

\begin{abstract}
This paper studies the distributed optimization problem in the presence of heavy-tailed gradient noises. Here, a heavy-tailed noise $\xi$ does not necessarily  adhere to the bounded variance assumption, i.e., $\textcolor{black}{\mathbb{E}}[\Vert \xi\Vert^2]\leq \nu^2$ for some positive constant $\nu$. Instead, it satisfies a more general condition  $\textcolor{black}{\mathbb{E}}[\Vert \xi\Vert^\delta]\leq \nu^\delta$ for some $1<\delta\leq 2$. The commonly-used bounded variance assumption is a special case of the considered noise assumption. While several distributed optimization algorithms have been proposed for scenarios involving heavy-tailed noise, these algorithms need a centralized server in the network which collects the information of all clients. Different from these algorithms, this paper considers a scenario where  there is no centralized server and the agents can only exchange information with neighbors within a communication graph. A distributed  method combining gradient clipping and distributed stochastic \textcolor{black}{gradient} projection is proposed. It is proven that when the gradient descent step-size and the gradient clipping step-size meet certain conditions, the state of each agent converges to an optimal solution of the distributed optimization problem with probability 1. The simulation results validate the algorithm.

\end{abstract}


\thispagestyle{empty}

\section{Introduction}

Distributed optimization is a method of achieving optimization through cooperation among multiple agents. The core idea of distributed optimization is to decompose a large problem into   small problems and distribute the local objectives to multiple intelligent agents for solution. Compared with centralized optimization, distributed optimization overcomes the disadvantages of single point of failure, high computation and communication burden, and restrictions on scalability and flexibility \cite{yang2019survey}. Distributed optimization has a wide range of applications, such as machine learning \cite{tsianos2012consensus}, distributed economic dispatch problem of energy systems\cite{yi2016initialization}, and formation control problem in robotic systems \cite{sun2021time}.

Different from deterministic optimization, stochastic optimization employs random factors to reach a solution to optimization problems. Stochastic optimization has \textcolor{black}{wide} applications in various fields. For example, the stochastic gradient descent (SGD) method  \textcolor{black}{can help escape from saddle points in nonconvex optimization and thus} has \textcolor{black}{significantly contributed to} the progress of \textcolor{black}{deep learning}. 

As an extension of distributed  deterministic  optimization to stochastic scenarios, distributed stochastic optimization has also attracted widespread attention. For instance, the authors in \cite{sundhar2010distributed} presented a distributed stochastic subgradient projection algorithm tailored for convex optimization challenges. \cite{srivastava2011distributed} introduced a distributed algorithm that relies on asynchronous step-sizes. \cite{wang2021distributed} proposed a proximal primal-dual approach geared towards nonconvex and nonsmooth problems. \cite{lu2023convergence} examined the high-probability convergence of distributed stochastic optimization.

Most of the above literature on stochastic optimization assume that the variance of the random \textcolor{black}{noise} is bounded. \textcolor{black}{Although} this is a mild assumption which covers many noises such as Gaussian noises, there are still many noises that do not satisfy this condition. For example, the Pareto distribution and the $\alpha$-stable \textcolor{black}{L$\acute{\text{e}}$vy} distribution \cite{zhang2020adaptive}. Heavy-tailed noises \textcolor{black}{have}  \textcolor{black}{recently been} observed in many machine learning systems \cite{zhang2020adaptive,yang2022taming,simsekli2019tail,gurbuzbalaban2021heavy,wang2021convergence,math13040665, gurbuzbalaban2022heavy,xu2023non}.  Moreover, the traditional SGD method may diverge when the noise is heavy-tailed \cite{zhang2020adaptive,yang2022taming}. \textcolor{black}{For distributed optimization problems, an interesting result was shown in \cite{gurbuzbalaban2022heavy} that the law of the distributed stochastic gradient decent iterates converges to a distribution with polynomially decaying heavy tails. However, \cite{gurbuzbalaban2022heavy} does not focus on how to tackle the heavy-tailed noise problem.}

In this work, we study the distributed stochastic optimization algorithm under the influence of heavy-tailed gradient noises using neighboring agents'  information exchange only. The main contributions of this paper are listed as follows: 

1) We design a distributed updating law to deal with the  stochastic optimization problem under heavy-tailed noises. Strict proof is given to show that under some mild parameter conditions, the state of each agent converges to an optimal solution with probability 1.

2) Compared with the existing distributed stochastic optimization works, \textcolor{black}{such as} \cite{sundhar2010distributed,srivastava2011distributed,alghunaim2019distributed, pu2021distributed,yuan2018optimal, wang2021distributed}, our algorithm allows \textcolor{black}{for} heavy-tailed noises, generalizing the commonly used assumption $\textcolor{black}{\mathbb{E}}[\Vert \xi\Vert^2]\leq \nu^2$. Specifically, our algorithm can handle noises that satisfy $\textcolor{black}{\mathbb{E}}[\Vert \xi\Vert^\delta]\leq \nu^\delta$ for some $1<\delta\leq 2$ and positive constant $\nu$. 

3) Compared with the existing literature on heavy-tailed noises, \textcolor{black}{such as}  \cite{zhang2020adaptive,gorbunov2020stochastic, liu2023high,cutkosky2021high, yang2022taming,  liu2022communication, yu2023smoothed, gorbunov2023high}, our algorithm is distributed in the sense that there is no central server and the agents can only exchange information via a strongly connected communication graph.  Note that despite that \cite{ yang2022taming, liu2022communication,yu2023smoothed, gorbunov2023high} studied heavy-tailed noises in a distributed setting, these algorithms require a central server \cite{ yang2022taming,yu2023smoothed, gorbunov2023high} or the agents require the states of all the other agents at some time instant \cite{liu2022communication}. These works are not related to multi-agent consensus while our work is related to consensus and graph theory.


\textbf{Notations}: Throughout this paper, $0$ represents a zero vector with an appropriate dimension or scalar zero.   $\mathbb{R}$ and \textcolor{black}{$\mathbb{R}^n$} represent the real number set and the \textcolor{black}{$n$}-dimensional real vector set, respectively. $\left\vert \cdot \right\vert$ is the absolute value of a scalar, and $\left\Vert \cdot \right\Vert$ is the 2-norm of a vector.  \textcolor{black}{``$\text{liminf}$" represents  the limit inferior of a sequence}. $\textcolor{black}{\mathbb{P}}_{\Omega}[\cdot]$ is the Euclidean projection of a vector to the space $\Omega$. $[e]_i$ is the $i$-th element of a vector $e$.  $[D]_{i,j}$ is the element in the $i$-th row and the $j$-th column of a matrix $D$. $\min\{a,b\}$ represents the minimum of $a$ and $b$. $\textcolor{black}{\mathbb{E}}[\cdot]$ is the expectation of a variable. \textcolor{black}{$\lfloor\cdot\rfloor$ is the rounding operation.} \textcolor{black} {A function $f(\cdot)$ is called $\mu$-strongly convex if and only if  $\forall x, y\in\mathbb{R}^n$, $f(y)-f(x)\geq \nabla f(x)(y-x)+\frac{\mu}{2}||y-x||^2,$ where $ \nabla f(x)$ is the gradient of $f(\cdot)$ at $x$ and $\mu>0$. }

%

\section{Problem Formulation\label{S3}}

Consider a multi-agent system comprised of  $N>1$ agents. The agents cooperate to solve the following problem 
\begin{flalign}
& \min_{\theta\in \Omega} f(\theta)=\sum_{i=1}^N f_i(\theta),& \label{system1}
\end{flalign}
where $\theta\in \Omega\subset\mathbb{R}^n$ is the decision variable, $f(\theta)$ is the global objective function,  $f_i(\theta)$ is the local objective function in agent $i$, and $\Omega$ is the local constraint set. 

Suppose that each agent can only get information from its neighbors via a communication graph $\mathcal{G}_{N}:=\{\mathcal{V}_{N},\mathcal{E}_{N}\}$ where $\mathcal{V}_{N}=\{1,\cdots,N\}$ is the node set, and    $\mathcal{E}_{N}\subset \mathcal{V}_{N}\times \mathcal{V}_{N}$ is the edge set.
\textcolor{black}{Let $A$ be a matrix with its element $[A]_{i,j}$ representing the weight of the edge from node $j$ to node $i$. The set $\mathcal{N}_i$ denotes the agents whose states are available to agent $i$ via one-hop communication. We assume that $i\in\mathcal{V}_i$. Furthermore, $j\in\mathcal{V}_i$ if and only if $[A]_{i,j}>0$.} 

\begin{assumption}\cite{sundhar2010distributed}
$\mathcal{G}_{N}$ is strongly connected, \textcolor{black}{i.e., there is a path between any two pair of vertices}. \textcolor{black} {The matrix $A$ is doubly stochastic, i.e., $\sum_{i=1}^N [A]_{i, j}= \sum_{j=1}^N [A]_{i, j}=1$}. There exists a scalar $0<\eta<1$ such that $[A]_{i, j}\geq \eta$ when $j\in\mathcal{N}_i$. \label{doubly}
\end{assumption}

\color{black}
\begin{assumption}
  $f_i(\theta)$ is continuously differentiable \textcolor{black}{and convex} \textcolor{black}{for $\theta\in\Omega$}.  \label{strongly}
\end{assumption}

\color{black}
\begin{assumption}
	The constraint set $\Omega$ is nonempty, closed and convex.  There exists a positive scalar $C_0$ such that $\Vert \nabla f_i(\theta) \Vert\leq C_0$ for all $\theta\in\Omega$.\label{boundedness}
\end{assumption}

\color{black}
\begin{assumption}
There exists at least one optimal solution $\theta^*\in\Omega^*$ to Problem \eqref{system1}, \textcolor{black}{where $\Omega^*$ is the optimal solution set}. \label{existence}
\end{assumption}

\begin{remark}
	Assumptions \ref{doubly},  \ref{strongly}, \ref{boundedness} and \ref{existence} are all commonly-used assumptions in the literature. For example, the doubly stochastic assumption is used in \cite{nedic2009distributed, sundhar2010distributed,yuan2018optimal,fallah2022robust}.
\end{remark}

\color{black}
\section{Main Results} \label{S4}
\textcolor{black}{In this section, we first design an algorithm based on consensus and gradient clipping techniques. Then, we prove the convergence of the proposed algorithm given some conditions. Finally, we provide the convergence rate of the algorithm.}
\subsection{Algorithm Design}

\textcolor{black}{Suppose that each agent has a local updating variable $	x_{i,k}$ at each iteration $k=0,1,...$ and $i\in\mathcal{V}_N.$ The objective is to design an updating law for agent $i$ using only neighboring states such that $	x_{i,k}$ converges to an optimal solution $\theta^*$.}

Motivated by \cite{sundhar2010distributed} and  \cite{zhang2020adaptive}, the distributed updating law for agent $i$ is designed as
\begin{flalign}
	&x_{i, k+1}=\textcolor{black}{\mathbb{P}}_\Omega\left[v_{i,k}-{\alpha_{k}}\hat{g}_{i,k}(v_{i,k})\right], x_{i, 0}\in\Omega \label{algorithm}&
\end{flalign}
where 	$v_{i,k}=\sum_{j=1}^N \textcolor{black}{[A]_{i, j}} x_{j, k}$, and 
\begin{flalign}
&\hat{g}_{i,k}(v_{i,k})=\min\left\{1,\frac{{\tau_{k}}}{\Vert g_{i,k}(v_{i,k})\Vert }\right\}g_{i,k}(v_{i,k}),&\notag\\
&g_{i,k}(v_{i,k})=\nabla f_i(v_{i,k})+\xi_{i, k}. &\label{al2}
\end{flalign}
In \eqref{algorithm}-\eqref{al2}, $\alpha_{k}$ and $\tau_{k}$ are positive sequences determined later, and $\xi_{i, k}$ is the random \textcolor{black}{noise} of agent $i$ at step $k$. 
\textcolor{black}{From \eqref{al2}, we have $\mathbb{E}[g_{i,k}(v_{i,k})]=\nabla f_i(v_{i,k})$. So Assumption \ref{boundedness} is equivalent to the bounded gradient condition.}

\color{black}
\begin{remark}
Note that there is another version of the distributed
algorithm, where the gradient is calculated on $x_{i,k}$ instead of
$v_{i,k}$  \cite{gurbuzbalaban2021decentralized}. We will explore the difference of the convergence rate between the two designs in future.
\end{remark}
\color{black}

\textcolor{black}{Let $\mathcal{F}_k$ be the \textcolor{black}{$\sigma$}-algebra $\sigma(\xi_{i, m};i\in\mathcal{V}_N, m=0,\cdots,k-1)$ generated
	by the errors up to $k-1$}. The following assumptions on the noises are made to facilitate the subsequent analysis.

\begin{assumption}  \cite{zhang2020adaptive} 
	There exist $\delta$ and  $\nu$ with $\delta\in(1,2]$ and $\nu>0$ such that  $\textcolor{black}{\mathbb{E}}[\Vert \xi_{i, k}\Vert^\delta|\mathcal{F}_k]\leq \nu^\delta$ with probability 1. \label{alpha}
\end{assumption}

%

\begin{assumption}
	(Unbiased Local Gradient Estimator)  $\textcolor{black}{\mathbb{E}}[\xi_{i, k}|\mathcal{F}_k]=0$ with probability 1. \label{unbiased}
\end{assumption}


%
\subsection{Convergence Analysis}

Let $b_{i,k}\textcolor{black}{:}=\hat{g}_{i,k}-\nabla f_i(v_{i,k})$ and  $B_{i,k}\textcolor{black}{:}=\textcolor{black}{\mathbb{E}}[b_{i,k}|\mathcal{F}_k]=\textcolor{black}{\mathbb{E}}[\hat{g}_{i,k}|\mathcal{F}_k]-\nabla f_i(v_{i,k})$. Based on the definition of $v_{i,k}$ and Assumption \ref{doubly}, $v_{i,k}\in \Omega $ for all $k\geq 0$. Thus, \textcolor{black}{by Assumptions \ref{strongly} and \ref{boundedness}}, $\Vert \nabla f_i(v_{i,k}) \Vert\leq C_0$ for all $k\geq 0$. 
%
%
Then, the following conclusion holds.
\begin{theorem}
	 Suppose that Assumptions \ref{doubly}-\ref{unbiased} hold, and \textcolor{black}{the decreasing sequence $\alpha_{k}$ and the increasing sequence} $\tau_{k}$ satisfy
	 \begin{flalign}
	 			& \textcolor{black} {\text{ there is a $k_1>0$ such that $\tau_{k}\geq 2C_0$ for all $k\geq k_1$,}} &\label{c00} \\
		 &\textcolor{black} {\sum_{k=0}^{+\infty}\alpha_{k}=\infty, \sum_{k=0}^{+\infty}\alpha_{k}^2<+\infty},  &\label{c1}\\
	 &\textcolor{black}{	\lim_{k\rightarrow +\infty}\alpha_{k}\tau_{k}=0},  \textcolor{black}{ \sum_{k=0}^{+\infty}\alpha_{k}^2 \tau_{k}^2< +\infty.}& \label{c3}
	 	\end{flalign}
Suppose that there exist constants $\omega>0$ and $\varpi\in\mathbb{R}$ such that
	 \begin{flalign}
	 	&\textcolor{black}{\sum_{k=0}^{+\infty} (\omega \alpha_k^\varpi-\mu\alpha_k)<+\infty}, 
	 		\textcolor{black}{	\sum_{k=0}^{+\infty} \alpha_k^{2-\varpi}\tau_k^{2-2\delta}<+\infty,}&\label{eq02}
	 	\end{flalign}
\textcolor{black}{where $\mu=0$ if $f_i(\theta)$, $\forall i \in\mathcal{V}_N$ is convex and $\mu>0$  if  $f_i(\theta)$, $\forall i \in\mathcal{V}_N$ is $\mu$-strongly convex.} Then, \textcolor{black}{under the algorithm given in \eqref{algorithm} and \eqref{al2}},
\begin{center}
\it{$x_{i, k}$ converges to a vector in $\Omega^*$ with probability 1},
\end{center} \textcolor{black}{where $\Omega^*$ is the optimal solution set to Problem \eqref{system1}}. \label{theo1}
\end{theorem}
\begin{proof}
\textcolor{black}{See Appendix \ref{app1}}. 
\end{proof}
\color{black}

\textcolor{black}{Moreover, by selecting specific stepsizes, we can obtain the following conclusions.}

\begin{proposition}
		Suppose that $f_i(\theta), \forall i\in\mathcal{V}_N$ is convex.  
	Let $\alpha_k=c(1+k)^{-a}$ and $\tau_k=d(1+k)^{b}$ with $c,d>0$. Then, any $a,b$ with $\frac{1}{2\delta}+\frac{1}{2}<a\leq 1$, $\frac{1}{a}<\varpi<\frac{2a-1}{a}\delta$,  $\frac{1+a\varpi-2a}{2\delta-2}<b<a-\frac{1}{2}$ satisfy the conditions in Theorem \ref{theo1}. \label{prop1}
	\end{proposition}

\begin{proof}
\textcolor{black}{See Appendix \ref{app2}}.
\end{proof}

\begin{proposition}
	Suppose that $f_i(\theta), \forall i\in\mathcal{V}_N$ is $\mu$-strongly convex. 	Let $\alpha_k=c(1+k)^{-a}$ and $\tau_k=d(1+k)^{b}$  with $c,d>0$. Then, any $a,b$ with $\frac{1}{2(2\delta-1)}+\frac{1}{2}<a\leq 1$, $1\leq \varpi<\frac{2a-1}{a}\delta$,  $\frac{1+a\varpi-2a}{2\delta-2}<b<a-\frac{1}{2}$ satisfy the conditions in Theorem \ref{theo1}. \label{prop2}
\end{proposition}
\begin{proof}
	\textcolor{black}{ See Appendix \ref{app3}}.
\end{proof}

\textcolor{black}{If the constraint set $\Omega$ is bounded, the conditions for the parameters can be relaxed, as shown below.}
\begin{theorem}
	Suppose that Assumptions \ref{doubly}-\ref{unbiased} hold and the constraint set $\Omega$ is bounded. Let \textcolor{black}{the decreasing sequence $\alpha_{k}$ and the increasing sequence $\tau_{k}$} satisfy conditions \eqref{c00}-\eqref{c3} and	
	\begin{flalign}
	&	\sum_{k=0}^{+\infty} \alpha_k\tau_k^{1-\delta}<+\infty. &\label{eq06}
	\end{flalign}
Then, \textcolor{black}{under the algorithm given in \eqref{algorithm} and \eqref{al2}}, 
\begin{center}
 \it{$x_{i, k}$ converges to a vector in $\Omega^*$ with probability 1}, 
 \end{center}
\textcolor{black}{where $\Omega^*$ is the optimal solution set to Problem \eqref{system1}}  \label{theorem2}
\end{theorem}
\begin{proof}
\textcolor{black}{See Appendix \ref{app4}}.
\end{proof}

\begin{proposition}
	Let $\alpha_k=c(1+k)^{-a}$ and $\tau_k=d(1+k)^{b}$ with $c,d>0$. Then,  any $a,b$ with $\textcolor{black}{\frac{1}{2\delta}+\frac{1}{2}}<a\leq 1$ and $\frac{1-a}{\delta-1}<b<a-\frac{1}{2}$ satisfy the conditions in Theorem \ref{theorem2}. \label{prop3}
\end{proposition}

\begin{proof}
\textcolor{black}{See Appendix \ref{app5}}.
\end{proof}

\subsection{Convergence Rate}
Motivated by \cite{wang2017distributed,xi2016distributed, pang2019randomized}, let $f_t=\min_{\lfloor\frac{t}{2}\rfloor\leq k\leq t}\textcolor{black}{\mathbb{E}}[f(y_{k})]$ where $\lfloor\frac{t}{2}\rfloor\geq \textcolor{black}{\max\{k_2, 3\}}$, $k_2 =\max\{2, k_1\}$, and $k_1$ was defined in \eqref{c00}. Then, we can arrive at the following conclusion.

\begin{theorem}
	Suppose that Assumptions \ref{doubly}-\ref{unbiased} hold. Let  $\alpha_k=c(1+k)^{-a}$ and $\tau_k=d(1+k)^{b}$ with $c,d>0$. (a) If  $f_i(\theta), \forall i\in\mathcal{V}_N$ is $\mu$-strongly convex,  $\frac{1}{2(2\delta-1)}+\frac{1}{2}<a<1$, and  $\frac{1-a}{2\delta-2}<b<a-\frac{1}{2}$, then 
$f_t-f(\theta^*)\leq \textcolor{black}{\frac{NC_2C_a}{2c(1+t)^{1-a}}}
\textcolor{black}{+\left(\frac{2C_0C_1C_aN^2\zeta}{1-\beta}+\frac{2cdC_0C_{\beta_1}C_aN^2\zeta}{(b-a+1)(1-\beta_1)}\right)\frac{1}{1+t}}
$ $\textcolor{black}{+\frac{4\beta cdC_0C_a N^2\zeta}{(2a-b-1)(1-\beta)(t-5)^{a-b}}}\textcolor{black}{+ \frac{d^{2-2\delta}NC_a(2\nu)^{2\delta}}{\mu(-b(2-2\delta)+a-1)(t-1)^{b(2\delta-2)}}}$ $\textcolor{black}{+\frac{4cdC_0C_aN}{(2a-b-1)(t-3)^{a-b}}}\textcolor{black}{+\frac{cd^2NC_a}{(2(a-b)-1)(t-1)^{a-2b}}
	=O(\frac{1}{t^{1-a}})}$. \newline (b) If the constraint set $\Omega$ is bounded,  $\textcolor{black}{\frac{1}{2\delta}+\frac{1}{2}}<a<1$ and $\frac{1-a}{\delta-1}<b<a-\frac{1}{2}$,  then 
		$f_t-f(\theta^*)\leq \textcolor{black}{\frac{NC_2C_a}{2c(1+t)^{1-a}}}
		\textcolor{black}{+\left(\frac{2C_0C_1C_aN^2\zeta}{1-\beta}+\frac{2cdC_0C_{\beta_1}C_aN^2\zeta}{(b-a+1)(1-\beta_1)}\right)\frac{1}{1+t}}
		$ $\textcolor{black}{+\frac{4 \beta cdC_0C_a N^2\zeta}{(2a-b-1)(1-\beta)(t-5)^{a-b}}}\textcolor{black}{+ \frac{d^{1-\delta}NC_2C_a(2\nu)^{\delta}}{(b(\delta-1)+a-1)(t-1)^{b(\delta-1)}}}$ $\textcolor{black}{+\frac{4cdC_0C_aN}{(2a-b-1)(t-3)^{a-b}}}\textcolor{black}{+\frac{cd^2NC_a}{(2(a-b)-1)(t-1)^{a-2b}}
			=O(\frac{1}{t^{1-a}})}$. \textcolor{black}{Here, $C_0, C_1, C_2$ are constants satisfying for all $i, k$, $\Vert \nabla f_i(v_{i,k}) \Vert\leq C_0$, $\Vert x_{i, 0}\Vert \leq C_1$, $\Vert v_{i,k}-\theta^*\Vert \leq C_2$, $C_a=\frac{1-a}{1-1.2^{a-1}}$,  $\sqrt{\beta}<\beta_1<1$, $C_{\beta_1}$ is a constant that is positively related to $\beta_1$, ${\zeta}=\left(1-\frac{\eta}{4N^2}\right)^{-2}$,  $\beta=\left(1-\frac{\eta}{4N^2}\right)^{\frac{1}{Q}}$, $\eta$ was defined in   Assumption \ref{doubly} and $Q$ is the number of edges in graph $\mathcal{G}_{N}$.} \label{conv1}
\end{theorem}

\begin{proof}
\textcolor{black}{See Appendix \ref{app6}}.
\end{proof}

\color{black}

\begin{remark}
Based on  \eqref{eqa2} and the definitions of $\beta$ and $\zeta$, \textcolor{black}{for the $\mu$-strongly convex case}, we can obtain that the upper bound of $f_t-f(\theta^*)$ has a positive correlation with the number of agents $N$, \textcolor{black}{the edge number $Q$ in the graph}, the bounded moment $\nu$ in Assumption \ref{alpha}, and has a negative correlation with the strong convex modulus $\mu$ and the tail index $\delta$ (for $\tau_{t-2}>2\nu$).
\end{remark}

%

%

\begin{remark}
	In \cite{pu2021sharp}, the convergence rate of $\textcolor{black}{\mathbb{E}}[||\bar{x}[k]-x^*||^2]$ with $\bar{x}[k]=\frac{1}{N}\sum_{i=1}^Nx_i[k]$ can be analyzed for a class of similar stochastic optimization algorithms with traditional noises. In future, 	inspired by \cite{pu2021sharp}, we will explore the convergence rate of $\textcolor{black}{\mathbb{E}}[||\bar{x}[k]-x^*||^2]$ for stochastic algorithms under heavy-tailed noises.
\end{remark}

\begin{remark}
	In this work, we assume that the gradient over the constraint set is bounded. In \cite{yuan2016convergence,pu2021sharp}, an analysis of the decentralized gradient decent algorithm without this assumption was given. In future, we will explore this aspect following this line.
\end{remark}

\color{black}
\begin{remark}
	When the objective function is convex, the optimal rate is approximately $O(t^{\frac{1}{2\delta}-\frac{1}{2}})$, which is obtained when $a$ is close to $\frac{1}{2\delta}+\frac{1}{2}$. When the objective function is $\mu$-strongly convex, the optimal rate is approximately $O(t^{\frac{1}{2(2\delta-1)}-\frac{1}{2}})$ when $a$ is close to $\frac{1}{2(2\delta-1)}+\frac{1}{2}$. Here, $\delta$ is the heavy-tail index. 
\end{remark}

\color{black}
\section{Simulation} \label{S5}

\textcolor{black}{In this section, we provide two examples to verify the effectiveness of the proposed algorithm.}

\subsection{Numerical Verification}

\color{black}
 In this section, we consider that $N=30$ agents collaborate to address a distributed optimization problem, where 
	$f_i(\theta)=\frac{1}{2}(\theta-0.01\times i\times\mathbf{1}_6)^\top(\theta-0.01\times i\times\mathbf{1}_6) $. The constraint set is given by $\Omega=\{\theta\in\mathbb{R}^6|\Vert\theta_i\Vert\leq 1, i=1,\cdots,6\}$. 
It can be verified that the optimization problem satisfies Assumptions \ref {strongly} and  \ref{boundedness}. The communication graph of the 30 agents is a circle graph. The initial values $x_i(0)=0$. Let \eqref{algorithm}-\eqref{al2} be the updating law with $\alpha_k=0.1{(k+1)^{-0.9}}$ and $\tau_k=5(k+1)^{0.3}$.
 \color{black}
 
 Let $\xi$ be a random \textcolor{black}{noise} with $\xi =\textcolor{black}{\phi} -2 \phi_{\min}$ where $\phi$ satisfies Pareto distribution with tail index $\gamma=2$ and $\phi_{\min}=1$.  The probability density function of $\phi$ can be described as $p(\phi)=\left\{\begin{array}{cc}0, & \text { if } \phi\leq \phi_{\min }, \\ \frac{2 \phi_{\min }^2}{\phi^{3}}, & \text { if } \phi >\phi_{\min } .\end{array}\right.$  The expectation of $\phi$ is $ {2 \phi_{\min}} $ and the variance is infinity. Thus, $\textcolor{black}{\mathbb{E}}[\xi]=\textcolor{black}{\mathbb{E}}[\phi]-2 \phi_{\min}=0$. Since $p(\xi)=\left\{\begin{array}{cc}0, & \text { if } \xi \leq - \phi_{\min} , \\ \frac{2 \phi_{\min }^2}{(\xi+ 2 \phi_{\min} )^{3}}, & \text { if } \xi >- \phi_{\min} ,\end{array}\right.$ we have for $\delta=\frac{3}{2}$, $\mathbb{E}[\vert\xi\vert^\delta]<3^\delta$.
 
Let  $\xi_{i,k}=\phi_{i,k} -2$ where each element of $\phi_{i,k}$ satisfies Pareto distribution with tail index 2 and the minimum parameter 1. Based on the above analysis, $\textcolor{black}{\mathbb{E}}[\xi_{i,k}]=0$ and  $\textcolor{black}{\mathbb{E}}[\Vert\xi_{i,k}\Vert^\delta]=\textcolor{black}{\mathbb{E}}[(\Vert\xi_{i,k}\Vert^2)^\frac{\delta}{2}]\leq \textcolor{black}{\mathbb{E}}[(6\max\vert{[\xi_{i,k}]_j}\vert^2)^\frac{\delta}{2}]<8^\delta $. Thus, $\xi_{i,k}$ satisfies Assumptions \ref{alpha} and \ref{unbiased}.

\color{black}
 
\begin{figure}[htb]
	\centering
	\includegraphics[width=6.8cm]{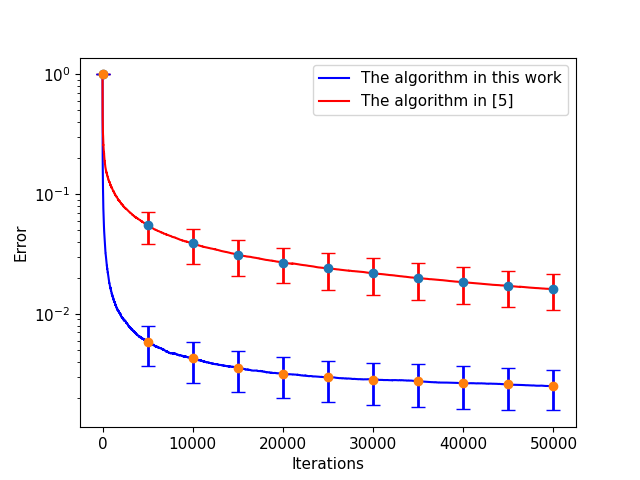} \caption{The evolution of $\text{log}_{10}\frac{f(y_k)-f(\theta^*)}{f(y_0)-f(\theta^*)}$.}
	\label{g2}
\end{figure}

Fig. \ref{g2} shows the evolution of $\text{log}_{10}\frac{f(y_k)-f(\theta^*)}{f(y_0)-f(\theta^*)}$ using the proposed algorithm and the distributed stochastic subgradient projection algorithm proposed in \cite{sundhar2010distributed}, where the vertical lines are the error bars. The optimal solution $\theta^*$ is found by the SLSQP method in the scipy.optimize package. Ten Monte Carlo trials are made and the results are averaged. As can be seen, the proposed algorithm has a smaller error.

\color{black}
 

%
%

\color{black}
\subsection{Real-World Dataset Simulation}

Inspired by \cite{gorbunov2020stochastic}, we consider the logistic regression problem for binary classification over real-world datasets Diabetes in the LibSVM library \cite{chang2011libsvm}. The gradient noise of the Diabetes dataset can be approximated by a heavy-tailed distribution \cite{gorbunov2020stochastic}. There are 4 agents in the circle graph. Each agent has a local database comprised of $192$ data points. The agent $i$'s objective function is 
$f_i(\theta)=\frac{1}{N_i}\sum_{\ell=1}^{N_i}\text{ln} (1+e^{(-a_\ell q_\ell^{\textcolor{black}{\top}}\theta)})$, 
where $q_\ell\in\mathbb{R}^8$ represents the feature of the $\ell$-th data point and $a_\ell\in\{-1,1\}$ represents the label of the $\ell$-th data point, and $N_i=10$ is the batch size.  The constraint set is given by $\Omega=\{\theta\in\mathbb{R}^8|\Vert\theta_i\Vert\leq 0.5, i=1,\cdots,8\}$.  The initial values $x_i(0)=0$. Let \eqref{algorithm}-\eqref{al2} be the updating law with $\alpha_k={3(k+1)^{-0.95}}$ and $\tau_k=5(k+1)^{0.4}$. 

\begin{figure}[htb]
	\centering
	\includegraphics[width=6.8cm]{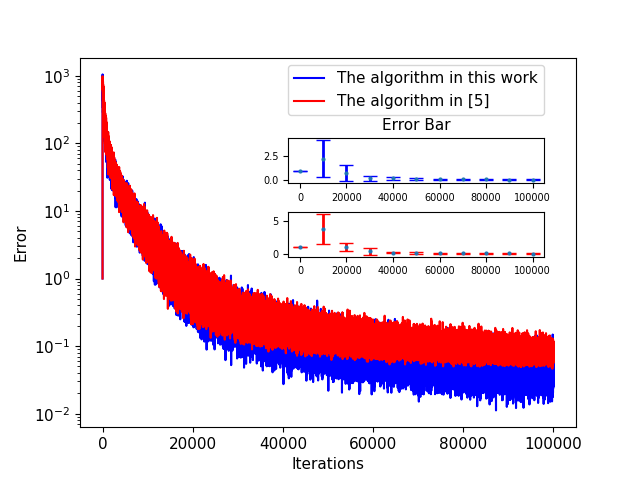} \caption{The evolution of $\text{log}_{10}\frac{f(y_k)-f(\theta^*)}{f(y_0)-f(\theta^*)}$.}
	\label{g3}
\end{figure}

Similarly, ten Monte Carlo trials are made and the results are averaged. Fig. \ref{g3} shows the evolution of $\text{log}_{10}\frac{f(y_k)-f(\theta^*)}{f(y_0)-f(\theta^*)}$ using the proposed algorithm and the distributed stochastic subgradient projection algorithm proposed in \cite{sundhar2010distributed}. It can be seen that the proposed algorithm has a better performance. 

\color{black}

\section{Conclusions} \label{S6}

In this paper, we considered the distributed optimization problem under the influence of heavy-tailed noises.  Unlike previous studies on distributed stochastic optimization, this paper considers the assumption of heavy-tailed distributions instead of the commonly considered  bounded variance noises, which is a broader assumption. \textcolor{black}{In future, we will extend the result to the nonconvex optimization scenario. }  \textcolor{black}{Specifically, we will consider the case where the objective function is nonconvex and $L$-smooth, and the case where the objective function is nonsmooth and weakly convex.}

\appendix
\subsection{Proof of Theorem \ref{theo1}} \label{app1}
\textcolor{black}{First, we define}  	$y_k\textcolor{black}{:}=\frac{1}{N} \sum_{i=1}^N x_{i, k}$ and
$p_{i, k+1}\textcolor{black}{:}=x_{i, k+1}-v_{i,k}=x_{i, k+1}-\sum_{j=1}^N \textcolor{black}{[A]_{i, j}}  x_{j, k}$. Then, we have
\begin{flalign}
	y_{k+1} & =\frac{1}{N}\sum_{i=1}^N\left( \sum_{j=1}^N \textcolor{black}{[A]_{i, j}}  x_{j, k}+ p_{i, k+1}\right) &\notag \\
	& =\frac{1}{N}\left(\sum_{j=1}^N\left(\sum_{i=1}^N \textcolor{black}{[A]_{i, j}} \right) x_{j, k}+\sum_{i=1}^N p_{i, k+1}\right)&\notag \\
	&= \frac{1}{N}\left(\sum_{j=1}^N x_{j, k}+\sum_{i=1}^N p_{i, k+1}\right)&\notag\\
	&=y_k+\frac{1}{N} \sum_{i=1}^N p_{i, k+1},& \notag
\end{flalign}
where \textcolor{black}{the second last equality is by} Assumption \ref{doubly}.
\textcolor{black}{Then, }
\begin{flalign}
	y_{k+1}&=y_0+\frac{1}{N} \sum_{\ell=1}^{k+1} \sum_{i=1}^N p_{i, \ell}&\notag\\
	&=\frac{1}{N} \sum_{i=1}^N x_{i, 0}+\frac{1}{N} \sum_{\ell=1}^{k+1} \sum_{i=1}^N p_{i, \ell}&\notag\\
	&=\frac{1}{N} \sum_{j=1}^N x_{j, 0}+\frac{1}{N} \sum_{\ell=1}^{k+1} \sum_{j=1}^N p_{j, \ell}.& \label{y}
\end{flalign}

Let $\mathbf{x}_k=[x_{1, k}^\top,\cdots,x_{N, k}^\top]^\top$, $\mathbf{p}_k=[p_{1, k}^\top,\cdots,p_{N, k}^\top]^\top$, $\tilde{A}=A\otimes I_n$. Then, 	$\mathbf{p}_{k+1}=\mathbf{x}_{k+1}-\tilde{A} \mathbf{x}_{k}$,
which implies that $\mathbf{x}_{k+1}=\mathbf{p}_{k+1}+\tilde{A} \mathbf{x}_{k} 
	=\mathbf{p}_{k+1}+\tilde{A}(\mathbf{p}_{k}+\tilde{A} \mathbf{x}_{k-1})=\sum_{\ell=0}^{k}\tilde{A}^\ell\mathbf{p}_{k+1-\ell}+\tilde{A}^{k+1}\mathbf{x}_0=\sum_{\ell=1}^{k+1}\tilde{A}^{k+1-\ell}\mathbf{p}_{\ell}+\tilde{A}^{k+1}\mathbf{x}_0,$
where in the last step we use $k+1-\ell$ to replace $\ell$, and $\tilde{A}^0:= I_{Nn}$. Thus, letting $A^0\triangleq I_{N}$ gives
\begin{flalign}
&	x_{i, k+1}
	=\sum_{\ell=1}^{k+1}\sum_{j=1}^{N}[A^{k+1-\ell}]_{i,j}p_{j, \ell}+\sum_{j=1}^{N}[A^{k+1}]_{i,j}x_{j,0}. &\label{x}
\end{flalign}
.

Based on \eqref{algorithm}, the definition of $p_{i,k+1}$ and the non-expansive property of the projection operator, we have
\begin{flalign}
	\Vert p_{i, k+1}\Vert& =\Vert \textcolor{black}{\mathbb{P}}_\Omega\left[v_{i,k}-\alpha_{k}\hat{g}_{i,k}\right]-v_{i,k}\Vert&\notag \\
	& \leq \alpha_{k}\Vert \hat{g}_{i,k}\Vert \leq \alpha_{k}\tau_{k}. &\label{p}
\end{flalign}
By \eqref{y} and \eqref{x}, we have for $k\geq 1$,
\begin{flalign}
	&\Vert y_{k+1}-x_{i, k+1}\Vert\notag\\
	=&\bigg\Vert \sum_{j=1}^N  \left(\frac{1}{N}-[A^{k+1}]_{i,j}\right)x_{j, 0} \notag \\ 
	&+\frac{1}{N} \sum_{\ell=1}^{k+1} \sum_{j=1}^N p_{j, \ell} - \sum_{\ell=1}^{k+1}\sum_{j=1}^N [A^{k+1-\ell}]_{i,j}p_{j, \ell}\bigg\Vert \notag \\
	=&\bigg\Vert \sum_{j=1}^N  \left(\frac{1}{N}-[A^{k+1}]_{i,j}\right)x_{j, 0}+\frac{1}{N} \sum_{\ell=1}^{k} \sum_{j=1}^N p_{j, \ell} \notag \\ 
	& - \sum_{\ell=1}^{k}\sum_{j=1}^N [A^{k+1-\ell}]_{i,j}p_{j, \ell}+ \frac{1}{N} \sum_{j=1}^N p_{j, k+1}-  p_{i, k+1}	\bigg\Vert \notag \\
	\leq& \sum_{j=1}^N  \left\vert\frac{1}{N}-[A^{k+1}]_{i,j}\right\vert \Vert x_{j, 0}\Vert \notag \\ 
	&+ \sum_{\ell=1}^{k} \sum_{j=1}^N \left\vert\frac{1}{N} -[A^{k+1-\ell}]_{i,j}\right\vert\Vert p_{j, \ell}\Vert \notag\\
	&+ \frac{1}{N} \sum_{j=1}^N \Vert p_{j, k+1}\Vert +\Vert p_{i, k+1} \Vert \notag\\
	\leq& N\zeta \beta^{k+1}C_1+ N \sum_{\ell=1}^{k}  \zeta \beta^{k+1-\ell} \alpha_{\ell-1}\tau_{\ell-1}+2\alpha_{k}\tau_{k},\label{y_x}
\end{flalign} 
\textcolor{black}{where in the last step we used Lemma 3.1 of \cite{sundhar2010distributed} and  \eqref{p}, ${\zeta}=\left(1-\frac{\eta}{4N^2}\right)^{-2}$,  $\beta=\left(1-\frac{\eta}{4N^2}\right)^{\frac{1}{Q}}$, $\eta$ was defined in   Assumption \ref{doubly}, $Q$ is the number of edges in graph $\mathcal{G}_{N}$} and $C_1$ satisfies $\Vert x_{j, 0}\Vert \leq C_1$ for $j\in\mathcal{V}_N$.
Then, for $k\geq 2$, we have
\begin{flalign}
	\Vert y_{k}-x_{i, k}\Vert\leq &
	N\zeta \beta^{k}C_1+ N \sum_{\ell=1}^{k-1}  \zeta \beta^{k-\ell} \alpha_{\ell-1}\tau_{ \ell-1}\notag\\
	&+2\alpha_{k-1}\tau_{ k-1}. &\label{y_x3}
\end{flalign}

For any $z\in\Omega$, based on the non-expansive property \textcolor{black}{of the projection operator} and the definition of $b_{i,k}$, \textcolor{black}{we have}
\begin{flalign}
	\Vert x_{i, k+1}-z\Vert ^2= & \Vert \textcolor{black}{\mathbb{P}}_\Omega\left[v_{i,k}-\alpha_{k}\hat{g}_{i,k}\right]-z\Vert^2 &\notag \\
	\leq & \Vert v_{i,k}-\alpha_{k}\hat{g}_{i,k}-z\Vert^2 &\notag \\
	= &  \Vert v_{i,k}-z\Vert^2-2 \alpha_{k}\hat{g}_{i,k}^\top \left(v_{i,k}-z\right)+\alpha_{k}^2\Vert\hat{g}_{i,k}\Vert^2 &\notag\\
	= & \Vert v_{i,k}-z\Vert^2-2 \alpha_{k}b_{i,k}^\top \left(v_{i,k}-z\right)&\notag \\
	& -2\alpha_{k} \nabla f_i(v_{i,k})^\top\left(v_{i,k}-z\right)+\alpha_{k}^2\Vert\hat{g}_{i,k}\Vert^2 &\notag\\
	\leq & (1-\mu\alpha_{k})\Vert v_{i,k}-z\Vert^2 &\notag\\
	&-2\alpha_{k} (f_i(v_{i,k})-f_i(z))&\notag \\
	&-2 \alpha_{k}b_{i,k}^\top \left(v_{i,k}-z\right)+\alpha_{k}^2\tau_{k}^2, &\label{eq1}
\end{flalign}
where in the last step we used the convexity of $f_i$ in Assumption \ref{strongly}. By the convexity of 2-norm and Assumption \ref{doubly}, 
\begin{flalign}
	\sum_{i=1}^N\Vert v_{i,k+1}-z\Vert^2&=\sum_{i=1}^N\bigg\Vert\sum_{j=1}^N \textcolor{black}{[A]_{i, j}}  x_{j, k+1}-z\bigg\Vert^2 &\notag\\
	&\leq \sum_{i=1}^N \sum_{j=1}^N \textcolor{black}{[A]_{i, j}} \Vert x_{j, k+1}-z\Vert ^2 &\notag\\
	&=\sum_{j=1}^N \sum_{i=1}^N \textcolor{black}{[A]_{i, j}} \Vert x_{j, k+1}-z\Vert ^2&\notag\\
	&=\sum_{j=1}^N\Vert x_{j, k+1}-z\Vert ^2.& \label{eq2}
\end{flalign}
Then, based on \eqref{eq1} and \eqref{eq2}, 
\begin{flalign}
	\sum_{i=1}^N\Vert v_{i,k+1}-z\Vert^2
	\leq & (1-\mu\alpha_{k} )\sum_{i=1}^N \Vert v_{i,k}-z\Vert^2 &\notag\\
	&-2\alpha_{k}\sum_{i=1}^N (f_i(v_{i,k})-f_i(z))&\notag \\
	&-2 \alpha_{k}\sum_{i=1}^Nb_{i,k}^\top \left(v_{i,k}-z\right)&\notag\\
	&+N\alpha_{k}^2\tau_{k}^2.& \label{vz}
\end{flalign}
Based on Assumption \ref{strongly}, we have
\begin{flalign}
	f_i\left(v_{i,k}\right)-f_i(z)  =&\left(f_i\left(v_{i,k}\right)-f_i\left(y_k\right)\right) \notag\\
	&+\left(f_i\left(y_k\right)-f_i(z)\right) & \notag\\
	\geq& -\Vert\nabla f_i\left(v_{i,k}\right)\Vert\Vert y_k-v_{i,k}\Vert \notag\\
	&+\left(f_i\left(y_k\right)-f_i(z)\right) &\notag\\
	=&  -\Vert\nabla f_i\left(v_{i,k}\right)\Vert\bigg\Vert y_k-\sum_{j=1}^N \textcolor{black}{[A]_{i, j}}  x_{j, k}\bigg\Vert \notag\\
	&+\left(f_i\left(y_k\right)-f_i(z)\right)  &\notag \\
	\geq&  -\Vert\nabla f_i\left(v_{i,k}\right)\Vert\sum_{j=1}^N \textcolor{black}{[A]_{i, j}} \Vert y_k- x_{j, k}\Vert &\notag\\
	&+\left(f_i\left(y_k\right)-f_i(z)\right).& \notag
\end{flalign}
Thus, according to Assumption \ref{doubly}, 
\begin{flalign}
&\hspace{-10pt}\sum_{i=1}^N	(f_i\left(v_{i,k}\right)-f_i(z)) &\notag\\
	\geq & -C_0\sum_{i=1}^N\sum_{j=1}^N \textcolor{black}{[A]_{i, j}} \Vert y_k- x_{j, k}\Vert+\sum_{i=1}^N\left(f_i\left(y_k\right)-f_i(z)\right)&\notag\\
	=&-C_0 \sum_{j=1}^N  \Vert y_k- x_{j, k}\Vert +\sum_{i=1}^N\left(f_i\left(y_k\right)-f_i(z)\right). &\notag
\end{flalign}
Furthermore, the inequality \eqref{vz} can be written as
\textit{\begin{flalign}
		\sum_{i=1}^N\Vert v_{i,k+1}-z\Vert^2
		\leq & (1-\mu\alpha_{k})\sum_{i=1}^N\Vert v_{i,k}-z\Vert^2 &\notag\\
		&+2C_0\alpha_{k} \sum_{j=1}^N  \Vert y_k- x_{j, k}\Vert&\notag \\
		&-2\alpha_{k}  \left(f \left(y_k\right)-f (z)\right)&\notag \\
		&-2 \alpha_{k}\sum_{i=1}^Nb_{i,k}^\top \left(v_{i,k}-z\right)&\notag\\
		&+N\alpha_{k}^2\tau_{k}^2.&\label{F1}
\end{flalign}}

\textcolor{black}{Taking the conditional expectations w.r.t. the $\sigma$-algebra $\mathcal{F}_k$ on both hand sides of} \eqref{F1} and letting $z=\theta^*$ gives
\begin{flalign}
	\sum_{i=1}^N\textcolor{black}{\mathbb{E}}[\Vert v_{i,k+1}-\theta^*\Vert^2| \mathcal{F}_k]
	\leq&(1-\mu\alpha_{k}) \sum_{i=1}^N\Vert v_{i,k}-\theta^*\Vert^2 &\notag\\
	&+2C_0\alpha_{k} \sum_{j=1}^N  \Vert y_k- x_{j, k}\Vert&\notag \\
	&-2\alpha_{k}  \left(f \left(y_k\right)-f (\theta^*)\right)&\notag \\
	&-2 \alpha_{k}\sum_{i=1}^NB_{i,k}^{\textcolor{black}{\top}} \left(v_{i,k}-\theta^*\right)&\notag\\
	&+N\alpha_{k}^2\tau_{k}^2.&\label{F2}
\end{flalign}


\textcolor{black}{According to Lemma C.1 of \cite{sadiev2023high}, for $k\geq k_1$ and for all $i$, $\Vert B_{i,k}\Vert \leq (2\nu)^\delta\tau_{k}^{1-\delta} $  with probability 1.} Then, utilizing Cauchy's inequality to the last but one term of \eqref{F2}, we have for $k\geq k_1$,
\begin{flalign}
	&	\sum_{i=1}^N\textcolor{black}{\mathbb{E}}[\Vert v_{i,k+1}-\theta^*\Vert^2| \mathcal{F}_k]&\notag\\
	\leq &    (\textcolor{black}{1+\omega \alpha_k^\varpi-\mu\alpha_k})\sum_{i=1}^N\Vert v_{i,k}-\theta^*\Vert^2 &\notag\\
	&+2C_0\alpha_{k} \sum_{j=1}^N  \Vert y_k- x_{j, k}\Vert-2\alpha_{k}  \left(f \left(y_k\right)-f (\theta^*)\right)&\notag \\
	&+\textcolor{black}{ \frac{N}{\omega}\alpha_k^{2-\varpi}(2\nu)^{2\delta} \tau_{k}^{2-2\delta}}+N\alpha_{k}^2\tau_k^2& \label{F3}
\end{flalign}
with probability 1.

According to \eqref{y_x3}, it can be obtained that for $k\geq 2$,
\begin{flalign}
	\sum_{k=2}^{+\infty}\alpha_{k} \sum_{j=1}^N  \Vert y_k- x_{j, k}\Vert 
	&\leq N^2\zeta C_1  \sum_{k=2}^{+\infty}	\alpha_{k} \beta^{k}&\notag \\
	&+ N^2   \zeta\sum_{k=2}^{+\infty}\alpha_{k}\sum_{\ell=1}^{k-1} \beta^{k-\ell} \alpha_{\ell-1}\tau_{\ell-1} &\notag\\
	&+2N\sum_{k=2}^{+\infty}\alpha_{k-1}^2\tau_{k-1}. &\label{y_x2}
\end{flalign}

Since $\sum_{k=1}^{+\infty}\alpha_k^2<+\infty$, $\alpha_k$ is bounded. According to the definition of $\beta$, $0<\beta<1$. Then, $N^2\zeta C_1  \sum_{k=2}^{+\infty}	\alpha_{k} \beta^{k}$ is bounded.

 Since $\sum_{k=0}^{+\infty}\alpha_{k}^2 \tau_{k}^2< +\infty $, based on Lemma 3.1 of \cite{sundhar2010distributed}, we have
	$N^2  \zeta\sum_{k=2}^{+\infty}\alpha_{k}\sum_{\ell=1}^{k-1} \beta^{k-\ell} \alpha_{\ell-1}\tau_{\ell-1}
	\leq 	\frac{N^2   \zeta}{2}\sum_{k=2}^{+\infty}\sum_{\ell=1}^{k-1}\beta^{k-\ell}(\alpha_k^2+\alpha_{\ell-1}^2\tau_{\ell-1}^2)
	\leq 	\frac{N^2   \zeta}{2}\sum_{k=2}^{+\infty}\alpha_k^2 \frac{\beta}{1-\beta}$ $+\sum_{k=2}^{+\infty}\sum_{\ell=1}^{k-1}\beta^{k-\ell}\alpha_{\ell-1}^2\tau_{\ell-1}^2
	<+\infty.$
	
Based on the conditions  $\sum_{k=0}^{+\infty}\alpha_{k}^2 \tau_{k}^2< +\infty$ and $\tau_{k}\geq 2C_0$ \textcolor{black}{for $k\geq k_1$}, we can get that $\sum_{k=2}^{+\infty}\alpha_{k-1}^2\tau_{k-1}<+\infty$ and $\sum_{k=2}^{+\infty}\alpha_{k}^2 \tau_{k} ^{2-\delta}<+\infty$. Thus, by \eqref{y_x2}, $	\sum_{k=2}^{+\infty}\alpha_{k} \sum_{j=1}^N  \Vert y_k- x_{j, k}\Vert <+\infty$.

Consequently, based on \eqref{F3}, condition  \eqref{c3} and \textcolor{black}{Theorem 1.3.12 of \cite{duflo2013random}}, $\Vert v_{i,k}-\theta^*\Vert$ converges and $\sum_{k=k_2}^{+\infty}\alpha_{k}  \left(f \left(y_k\right)-f (\theta^*)\right)<+\infty$ with probability 1, where $k_2=\max\{2,k_1\}$. Based on \eqref{y_x}, the condition  in the theorem and Lemma 3.1 of \cite{sundhar2010distributed}, $\lim_{k\rightarrow+\infty}\Vert y_k- x_{i, k}\Vert=0$. Thus, $\lim_{k\rightarrow+\infty}\Vert   y_k-v_{i,k}\Vert=0$, which implies that $\Vert y_k-\theta^*\Vert$ converges with probability 1. Since  $\sum_{k=0}^{+\infty}\alpha_k=\infty$ and  $\sum_{k=k_2}^{+\infty}\alpha_{k}  \left(f \left(y_k\right)-f (\theta^*)\right)<+\infty$ with probability 1, we have $\textcolor{black}{\liminf_{k\rightarrow+\infty}}f( y_k)=f (\theta^*)$ with probability 1. Based on the continuity of $f$, $y_k$  converges to  $\theta^*$ with probability 1. Then, $x_{i, k}$ converges to  $\theta^*$ with probability 1.
\subsection{Proof of Proposition \ref{prop1}}\label{app2}
 \textcolor{black}{Since $\frac{1}{2}<a\leq 1$, \eqref{c1} is satisfied. Since $b<a-\frac{1}{2}$,  \eqref{c3} is satisfied. Since $\varpi>\frac{1}{a}$ and $b>\frac{1+a\varpi-2a}{2\delta-2}$,  \eqref{eq02} is satisfied. The condition $\frac{1}{a}<\varpi<\frac{2a-1}{a}\delta$ can guarantee the existence of $b$, i.e., $\frac{1+a\varpi-2a}{2\delta-2}<a-\frac{1}{2}$. The condition $a>\frac{1}{2\delta}+\frac{1}{2}$ can guarantee the existence of $\varpi$, i.e., $\frac{1}{a}<\frac{2a-1}{a}\delta$. }
 
 \subsection{Proof of Proposition \ref{prop2}}\label{app3}
 	 \textcolor{black}{Similar to the convex case, the conditions $\frac{1}{2}<a\leq 1$,  $b<a-\frac{1}{2}$, $\varpi\geq 1$, $b>\frac{1+a\varpi-2a}{2\delta-2}$  can guarantee that \eqref{c1}--\eqref{eq02} hold. The condition $1\leq \varpi<\frac{2a-1}{a}\delta$ can  guarantee the existence of $b$, i.e., $\frac{1+a\varpi-2a}{2\delta-2}<a-\frac{1}{2}$.  The condition $a>\frac{1}{2(2\delta-1)}+\frac{1}{2}$ can  guarantee the existence of $\varpi$, i.e., $1<\frac{2a-1}{a}\delta$. }
\subsection{Proof of Theorem \ref{theorem2}}\label{app4}
 	 	From \eqref{F2} and the boundedness of $(v_{i,k}-\theta^*)$, we can obtain that
 	 for $k\geq k_1$,
 	 \begin{flalign}
 	 		\sum_{i=1}^N&\textcolor{black}{\mathbb{E}}[\Vert v_{i,k+1}-\theta^*\Vert^2| \mathcal{F}_k]&\notag\\
 	 	\leq & (1-\mu\alpha_k)   \sum_{i=1}^N\Vert v_{i,k}-\theta^*\Vert^2 +2C_0\alpha_{k} \sum_{j=1}^N  \Vert y_k- x_{j, k}\Vert&\notag \\
 	 	&-2\alpha_{k}  \left(f \left(y_k\right)-f (\theta^*)\right)
 	 	+\textcolor{black}{ NC_2\alpha_k(2\nu)^{\delta} \tau_{k}^{1-\delta}} &  \notag\\
 	 	&+N\alpha_{k}^2\tau_k^2 \label{F6}
 	 \end{flalign}
 	 with probability 1 where $C_2$ is a constant satisfying $\Vert v_{i,k}-\theta^*\Vert \leq C_2$ for all $v_{i,k}\in\Omega$. The following proof is similar to that of Theorem \ref{theo1} and thus is omitted.
\subsection{Proof of Proposition \ref{prop3}}\label{app5} 
 	 	\textcolor{black}{The conditions $\frac{1}{2}<a\leq 1$,   $b<a-\frac{1}{2}$, $-a+b(1-\delta)<-1$ can guarantee that \eqref{c1}, \eqref{c3}, and \eqref{eq06} hold, respectively. The condition $a>\textcolor{black}{\frac{1}{2\delta}+\frac{1}{2}}$ can guarantee the existence of $b$. }
 	 	
 \subsection{Proof of Theorem \ref{conv1}}\label{app6} 	 	
 	 	
 	 		\textcolor{black}{First, we prove case (a). } According to \eqref{y_x3}, selecting $\varpi=1$ and $\omega=\mu$, and taking the total expectation of both sides of \eqref{F3}, we have for $k\geq k_2$,
 	 		$\sum_{i=1}^N\textcolor{black}{\mathbb{E}}[\Vert v_{i,k+1}-\theta^*\Vert^2]
 	 		\leq\sum_{i=1}^N\textcolor{black}{\mathbb{E}}[\Vert v_{i,k}-\theta^*\Vert^2] +E_{0,k}
 	 		-2\alpha_{k} \textcolor{black}{\mathbb{E}}[\left(f \left(y_k\right)-f (\theta^*)\right)],$
 	 	where $E_{0,k}$ is defined by
 	 	\begin{flalign}
 	 		E_{0,k}=&2C_0\alpha_{k}N(	N\zeta \beta^{k}C_1+ N \sum_{\ell=1}^{k-1}  \zeta \beta^{k-\ell} \alpha_{\ell-1}\tau_{ \ell-1}+2\alpha_{k-1}&\notag\\
 	 		&\cdot\tau_{ k-1})
 	 		+\frac{N}{\mu}\alpha_k(2\nu)^{2\delta}\tau_{k}^{2-2\delta}+N\alpha_{k}^2\tau_k^2. &\notag
 	 	\end{flalign}
 	 		Summing from $k=\lfloor \frac{t}{2}\rfloor$ to $t$ gives
 	 	\begin{flalign}
 	 		 f_t-f(\theta^*) 
 	 		\leq&\frac{\sum_{i=1}^N\textcolor{black}{\mathbb{E}}[\Vert v_{i,\lfloor \frac{t}{2}\rfloor}-\theta^*\Vert^2]+ \sum_{k=\lfloor \frac{t}{2}\rfloor}^{t}E_{0,k}}{	\sum_{k=\lfloor \frac{t}{2}\rfloor}^{t}2\alpha_{k} }. & \label{eqa2}
 	 	\end{flalign}
 	 	
 	 	Next, we analyze the convergence rate term by term. It can be obtained that
 	 	\begin{flalign}
 	 		\frac{\sum_{k=\lfloor\frac{t}{2}\rfloor}^t\alpha_k\beta^k}{\sum_{k=\lfloor \frac{t}{2}\rfloor}^t\alpha_k}	=&
 	 		\frac{\sum_{k=\lfloor \frac{t}{2}\rfloor}^t(1+k)^{-a}\beta^k}{\sum_{k=\lfloor \frac{t}{2}\rfloor}^t(1+k)^{-a}} &\notag\\
 	 		\leq&  \frac{(1+\lfloor \frac{t}{2}\rfloor)^{-a}\frac{1}{1-\beta}}{\int_{\lfloor \frac{t}{2}\rfloor}^t(1+\ell)^{-a}d\ell}&\notag\\
 	 		\leq&  \frac{(1+\lfloor \frac{t}{2}\rfloor)^{-a}\frac{1}{1-\beta}}{\frac{(1+t)^{1-a}-(1+\lfloor \frac{t}{2}\rfloor)^{1-a}}{1-a}}&\notag\\
 	 		\leq & \frac{(\frac{1+t}{2})^{-a}\frac{1}{1-\beta}}{\frac{(1+t)^{1-a}-(\frac{1+t}{1.2})^{1-a}}{1-a}} &\notag \\
 	 		\leq &\textcolor{black}{ \frac{2(1-a)}{(1-{1.2^{a-1}})(1-\beta)(1+t)}}.& \label{bounded2}
 	 	\end{flalign}
 	 	Similarly, we can get that
 	 	\begin{flalign}
 	 		\frac{\sum_{i=1}^N\textcolor{black}{\mathbb{E}}[\Vert v_{i,\lfloor \frac{t}{2}\rfloor}-\theta^*\Vert^2]}{	\sum_{k=\lfloor \frac{t}{2}\rfloor}^{t}\alpha_{k} }
 	 		\leq &\frac{\sum_{i=1}^N\textcolor{black}{\mathbb{E}}[\Vert v_{i,\lfloor \frac{t}{2}\rfloor}-\theta^*\Vert^2]}{\textcolor{black}{c}\frac{(1+t)^{1-a}-(\frac{1+t}{1.2})^{1-a}}{1-a}} &\notag\\
 	 		\leq & \textcolor{black}{\frac{NC_2(1-a)}{c(1-{1.2^{a-1}})(1+t)^{1-a}}}. \label{bounded1}
 	 	\end{flalign}
 	 	Moreover,
 	 	\begin{flalign}
 	 	&\hspace{-10pt}\sum_{\ell=1}^{k-1}  \beta^{k-\ell} \alpha_{\ell-1}\tau_{ \ell-1}	=cd\sum_{\ell=1}^{k-1}  \beta^{k-\ell} \ell^{b-a}&\notag\\
 	 		=&\beta cd\sum_{\ell=2}^{\lfloor \frac{k}{2}\rfloor}  \beta^{k-\ell} (\ell-1)^{b-a}+\beta cd\sum_{\ell=\lfloor \frac{k}{2}\rfloor+1}^{k}  \beta^{k-\ell} (\ell-1)^{b-a}.& \label{eq1add}
 	 	\end{flalign}
 	 	
 	 	\color{black}
 	 	According to the conditions in the theorem, $\frac{1}{2}<a-b<1$, $0<b-a+1<\frac{1}{2}$ and $-1<b-2a+1<0$. For the first term in \eqref{eq1add}, motivated by \cite{wang2017distributed}, we have
 	 	$	\beta cd\sum_{\ell=2}^{\lfloor \frac{k}{2}\rfloor}  \beta^{k-\ell} (\ell-1)^{b-a} 
 	 		\leq \beta cd\beta^{k-\lfloor \frac{k}{2}\rfloor}\sum_{\ell=2}^{\lfloor \frac{k}{2}\rfloor} (\ell-1)^{b-a}
 	 		\leq \beta cd\beta^{\frac{k}{2}}\int_{\ell=1}^{\lfloor \frac{k}{2}\rfloor}(\ell-1)^{b-a}d\ell
 	 		= \frac{\beta cd\beta^{\frac{k}{2}}}{b-a+1}(\lfloor \frac{k}{2}\rfloor-1)^{b-a+1}.$
 	 	For any $\sqrt{\beta}<\beta_1<1$, $\lim_{k\rightarrow +\infty}\frac{\beta^{\frac{k}{2}+1}((\lfloor \frac{k}{2}\rfloor-1)^{b-a+1})}{\beta_1^k}=0$, which implies that there exists a constant $C_{\beta_1}$ such that $\beta cd\sum_{\ell=2}^{\lfloor \frac{k}{2}\rfloor}  \beta^{k-\ell} (\ell-1)^{b-a}\leq \frac{cdC_{\beta_1}\beta_1^k}{b-a+1}$.
 	 	Moreover, $
 	 		\beta cd\sum_{\ell=\lfloor \frac{k}{2}\rfloor+1}^{k}  \beta^{k-\ell} (\ell-1)^{b-a}
 	 	\leq\beta cd \lfloor \frac{k}{2}\rfloor^{b-a}\sum_{\ell=\lfloor \frac{k}{2}\rfloor+1}^{k}  \beta^{k-\ell}
 	 		\leq\frac{\beta cd}{1-\beta} (\frac{k-1}{2})^{b-a}.$ 
 	 		\color{black}
 	 Thus, we can obtain that
 	 \color{black}
 	 	\begin{flalign}
 	 		&\hspace{-10pt}\frac{\sum_{k=\lfloor\frac{t}{2}\rfloor}^t\alpha_k\sum_{\ell=1}^{k-1}  \beta^{k-\ell} \alpha_{\ell-1}\tau_{ \ell-1}}{\sum_{k=\lfloor\frac{t}{2}\rfloor}^t\alpha_k} &\notag\\
 	 		\leq 	&
 	 		\frac{\sum_{k=\lfloor\frac{t}{2}\rfloor}^t(1+k)^{-a}(\frac{cdC_{\beta_1}\beta_1^k}{b-a+1}+\frac{\beta cd}{1-\beta} (\frac{k-1}{2})^{b-a})}{\sum_{k=\lfloor\frac{t}{2}\rfloor}^t(1+k)^{-a}} &\notag\\
 	 		\leq &\textcolor{black}{ \frac{2cdC_{\beta_1}(1-a)}{(b-a+1)(1-{1.2^{a-1}})(1-\beta_1)(1+t)}}&\notag\\
 	 		&+ \frac{\frac{\beta cd}{1-\beta} (\frac{1}{2})^{b-a}\int_{\lfloor\frac{t}{2}\rfloor-1}^{t}(\ell-1)^{b-2a}d\ell}{\frac{(1+t)^{1-a}-(\frac{1+t}{1.2})^{1-a}}{1-a}}&\notag\\
 	 		\leq &\textcolor{black}{ \frac{2cdC_{\beta_1}(1-a)}{(b-a+1)(1-{1.2^{a-1}})(1-\beta_1)(1+t)}} &\notag\\
 	 		&+\frac{\frac{\beta cd}{1-\beta} (\frac{1}{2})^{b-a}\frac{1}{2a-b-1}(\lfloor\frac{t}{2}\rfloor-2)^{b-2a+1})}{\frac{(1+t)^{1-a}-(\frac{1+t}{1.2})^{1-a}}{1-a}}&\notag\\
 	 		\leq& \textcolor{black}{ \frac{2cdC_{\beta_1}(1-a)}{(b-a+1)(1-{1.2^{a-1}})(1-\beta_1)(1+t)}}&\notag\\
 	 		&\textcolor{black}{+\frac{4\beta cd(1-a)}{(2a-b-1)(1-{1.2^{a-1}}){(1-\beta)}(t-5)^{a-b}}.} &\label{bounded3}
 	 	\end{flalign}
 	 	\color{black}
 	 	
 	 	It can be obtained that
 	 	\begin{flalign}
 	 		&\hspace{-10pt}\frac{\sum_{k=\lfloor\frac{t}{2}\rfloor}^t\alpha_k\alpha_{k-1}\tau_{k-1}}{\sum_{k=\lfloor \frac{t}{2}\rfloor}^t\alpha_k}=
 	 		\frac{\textcolor{black}{cd}\sum_{k=\lfloor \frac{t}{2}\rfloor}^t(1+k)^{-a}k^{-a}k^{b}}{\sum_{k=\lfloor \frac{t}{2}\rfloor}^t(1+k)^{-a}} &\notag\\
 	 		\leq& \frac{\textcolor{black}{cd}\sum_{k=\lfloor \frac{t}{2}\rfloor}^tk^{b-2a}}{\sum_{k=\lfloor \frac{t}{2}\rfloor}^t(1+k)^{-a}} &\notag\\
 	 		\leq & \textcolor{black}{\frac{2cd(1-a)}{(2a-b-1)(1-{1.2^{a-1}})(t-3)^{a-b}}.} &\label{bounded4}
 	 	\end{flalign}
 	 	Since $b(2-2\delta)-a<-1$, we have
 	 	\begin{flalign}
 	 		&\hspace{-10pt}\frac{\sum_{k=\lfloor\frac{t}{2}\rfloor}^t\alpha_k\tau_{k}^{2-2\delta}}{\sum_{k=\lfloor \frac{t}{2}\rfloor}^t\alpha_k}=
 	 		\frac{\textcolor{black}{d^{2-2\delta}}\sum_{k=\lfloor \frac{t}{2}\rfloor}^t(1+k)^{-a}(1+k)^{b(2-2\delta)}}{\sum_{k=\lfloor \frac{t}{2}\rfloor}^t(1+k)^{-a}} &\notag\\
 	 	\leq& \frac{\frac{-\textcolor{black}{d^{2-2\delta}}}{b(2-2\delta)-a+1} (\lfloor \frac{t}{2}\rfloor)^{b(2-2\delta)-a+1}}{\frac{(1+t)^{1-a}-(\frac{1+t}{1.2})^{1-a}}{1-a}}&\notag\\
 	 		\leq& \textcolor{black}{\frac{2d^{2-2\delta}(1-a)}{(-b(2-2\delta)+a-1)(1-{1.2^{a-1}})(t-1)^{b(2\delta-2)}}}. \notag
 	 	\end{flalign}
 	 	Since $1<2(a-b)<2$ and $0<a-2b<1$, we have
 	 	\begin{flalign}
 	 		&\hspace{-10pt}\frac{\sum_{k=\lfloor\frac{t}{2}\rfloor}^t\alpha_k^2\tau_{k}^{2}}{\sum_{k=\lfloor \frac{t}{2}\rfloor}^t\alpha_k}=
 	 		\frac{cd^2\sum_{k=\lfloor \frac{t}{2}\rfloor}^t(1+k)^{-2(a-b)}}{\sum_{k=\lfloor \frac{t}{2}\rfloor}^t(1+k)^{-a}} &\notag\\
 	 		\leq& \frac{\frac{-cd^2}{1-2(a-b)} (\lfloor \frac{t}{2}\rfloor)^{1-2(a-b)}}{\frac{(1+t)^{1-a}-(\frac{1+t}{1.2})^{1-a}}{1-a}}&\notag\\
 	 		\leq &\textcolor{black}{\frac{2cd^2(1-a)}{(2(a-b)-1)(1-{1.2^{a-1}})(t-1)^{a-2b}}}. &\label{pra}
 	 	\end{flalign}
 	 	Thus, based on \eqref{eqa2}, we can get that 
 	 	\begin{flalign}
 	 		&\hspace{-10pt}f_t-f(\theta^*)&\notag\\ \leq&\textcolor{black}{\frac{NC_2C_a}{2c(1+t)^{1-a}}+\frac{2C_0C_1C_aN^2\zeta}{(1-\beta)(1+t)}}&\notag\\
 	 	&\textcolor{black}{+\frac{2cdC_0C_{\beta_1}C_aN^2\zeta}{(b-a+1)(1-\beta_1)(1+t)}}&\notag\\
 	 	&\textcolor{black}{+\frac{4\beta cdC_0C_a N^2\zeta}{(2a-b-1)(1-\beta)(t-5)^{a-b}}}&\notag\\
 	 	&\textcolor{black}{+ \frac{d^{2-2\delta}NC_a(2\nu)^{2\delta}}{\mu(-b(2-2\delta)+a-1)(t-1)^{b(2\delta-2)}}}&\notag\\ &\textcolor{black}{+\frac{4cdC_0C_aN}{(2a-b-1)(t-3)^{a-b}}}&\notag\\
 	 	&\textcolor{black}{+\frac{cd^2NC_a}{(2(a-b)-1)(t-1)^{a-2b}}}.&\label{29}
 	 \end{flalign}
 	 	
 	 	\color{black}
 	 	\textcolor{black}{The proof for case (b) is similar to that of case (a).} According to \eqref{F6} and the boundedness of $\Omega$, we have	$\sum_{i=1}^N {\mathbb{E}}[\Vert v_{i,k+1}-\theta^*\Vert^2]
 	 		\leq\sum_{i=1}^N {\mathbb{E}}[\Vert v_{i,k}-\theta^*\Vert^2] +E_{1,k}
 	 		-2\alpha_{k}  {\mathbb{E}}[\left(f \left(y_k\right)-f (\theta^*)\right)],$
 	 	where $E_{1,k}=2C_0\alpha_{k}N(	N\zeta \beta^{k}C_1+ N \sum_{\ell=1}^{k-1}  \zeta \beta^{k-\ell} \alpha_{\ell-1}\tau_{ \ell-1}+2\alpha_{k-1}\tau_{ k-1})+ NC_2\alpha_k(2\nu)^{\delta} \tau_{k}^{1-\delta}+N\alpha_{k}^2\tau_k^2. $
 	 		Since $b(1-\delta)-a<-1$, we have
 	 	\begin{flalign}
 	 		\frac{\sum_{k=\lfloor\frac{t}{2}\rfloor}^t\alpha_k\tau_{k}^{1-\delta}}{\sum_{k=\lfloor \frac{t}{2}\rfloor}^t\alpha_k}=&
 	 		\frac{d^{1-\delta}\sum_{k=\lfloor \frac{t}{2}\rfloor}^t(1+k)^{-a}(1+k)^{b(1-\delta)}}{\sum_{k=\lfloor \frac{t}{2}\rfloor}^t(1+k)^{-a}} &\notag\\
 	 		\leq& \frac{\frac{-d^{1-\delta}}{b(1-\delta)-a+1} (\lfloor \frac{t}{2}\rfloor)^{b(1-\delta)-a+1}}{\frac{(1+t)^{1-a}-(\frac{1+t}{1.2})^{1-a}}{1-a}}&\notag\\
 	 		\leq&\frac{2d^{1-\delta}C_a}{(b(\delta-1)+a-1)(t-1)^{b(\delta-1)}}.&\label{last}
 	 	\end{flalign}
 	 	Then, similar to \eqref{eqa2} and \eqref{29}, we can get that
 	 	\begin{flalign}
	&\hspace{-10pt}f_t-f(\theta^*)&\notag\\ \leq&\textcolor{black}{\frac{NC_2C_a}{2c(1+t)^{1-a}}+\frac{2C_0C_1C_aN^2\zeta}{(1-\beta)(1+t)}}&\notag\\
&\textcolor{black}{+\frac{2cdC_0C_{\beta_1}C_aN^2\zeta}{(b-a+1)(1-\beta_1)(1+t)}}&\notag\\
&\textcolor{black}{+\frac{4\beta cdC_0C_a N^2\zeta}{(2a-b-1)(1-\beta)(t-5)^{a-b}}}&\notag\\
&\textcolor{black}{+ \frac{d^{1-\delta}NC_2C_a(2\nu)^{\delta}}{(b(\delta-1)+a-1)(t-1)^{b(\delta-1)}}}&\notag\\ &\textcolor{black}{+\frac{4cdC_0C_aN}{(2a-b-1)(t-3)^{a-b}}}&\notag\\
&\textcolor{black}{+\frac{cd^2NC_a}{(2(a-b)-1)(t-1)^{a-2b}}},\notag
		\end{flalign}
 	 	where we use  \eqref{bounded2}, \eqref{bounded1}, \eqref{bounded3}, \eqref{bounded4}, \eqref{pra} and \eqref{last}.

 \color{black}

\bibliographystyle{IEEEtran}
\bibliography{bib}
\vspace{1cm}
\end{document}